%
%

\def\juerg #1{#1}

\def\revis#1{#1}



\def\LaTeX{%
  \let\Begin\begin
  \let\End\end
  \def\Bcenter{\Begin{center}}
  \def\Ecenter{\End{center}}
  \let\Label\label
  \let\salta\relax
  \let\finqui\relax
  \let\futuro\relax}

\def\UK{\def\our{our}\let\sz s}
\def\USA{\def\our{or}\let\sz z}



\LaTeX

\USA


\salta

\documentclass[twoside,12pt]{article}
\setlength{\textheight}{24cm}
\setlength{\textwidth}{16cm}
\setlength{\oddsidemargin}{2mm}
\setlength{\evensidemargin}{2mm}
\setlength{\topmargin}{-15mm}
\parskip2mm


%
%
\usepackage{cite}

\usepackage{color}
\usepackage{amsmath}
\usepackage{amsthm}
\usepackage{amssymb}

\usepackage{amsfonts}
\usepackage{mathrsfs}

\usepackage{hyperref}
\usepackage[mathcal]{euscript}

\usepackage[ulem=normalem,draft]{changes}




%

\finqui

\def\Beq{\Begin{equation}}
\def\Eeq{\End{equation}}

\def\Bthm{\Begin{theorem}}
\def\Ethm{\End{theorem}}

\def\Brem{\Begin{remark}\rm}
\def\Erem{\End{remark}}

\def\Bdim{\Begin{proof}}
\def\Edim{\End{proof}}
\let\non\nonumber




\def\step #1 \par{\medskip\noindent{\bf #1.}\quad}


\def\Lip{Lip\-schitz}

\def\lhs{left-hand side}
\def\rhs{right-hand side}



\def\multibold #1{\def\arg{#1}%
  \ifx\arg\pto \let\next\relax
  \else
  \def\next{\expandafter
    \def\csname #1#1#1\endcsname{{\bf #1}}%
    \multibold}%
  \fi \next}

\def\pto{.}

\def\multical #1{\def\arg{#1}%
  \ifx\arg\pto \let\next\relax
  \else
  \def\next{\expandafter
    \def\csname cal#1\endcsname{{\cal #1}}%
    \multical}%
  \fi \next}


\def\multimathop #1 {\def\arg{#1}%
  \ifx\arg\pto \let\next\relax
  \else
  \def\next{\expandafter
    \def\csname #1\endcsname{\mathop{\rm #1}\nolimits}%
    \multimathop}%
  \fi \next}

\multibold
qwertyuiopasdfghjklzxcvbnmQWERTYUIOPASDFGHJKLZXCVBNM.

\multical
QWERTYUIOPASDFGHJKLZXCVBNM.

\multimathop
dist div dom meas sign supp .


\def\Accorpa #1#2 #3 {\gdef #1{\eqref{#2}--\eqref{#3}}%
  \wlog{}\wlog{\string #1 -> #2 - #3}\wlog{}}


\def\<#1>{\mathopen\langle #1\mathclose\rangle}
\def\norma #1{\mathopen \| #1\mathclose \|}


\def\iO{\int_\Omega}

\def\dt{\partial_t}
\def\dtt{\partial_{tt}}
\def\dn{\partial_{\bf n}}

\def\cpto{\,\cdot\,}

\def\checkmmode #1{\relax\ifmmode\hbox{#1}\else{#1}\fi}

\def\aeQ{\checkmmode{a.e.\ in~$Q$}}

\def\aat{\checkmmode{for a.e.~$t\in(0,T)$}}


\def\erre{{\mathbb{R}}}
\def\enne{{\mathbb{N}}}




\def\genspazio #1#2#3#4#5{#1^{#2}(#5,#4;#3)}
\def\spazio #1#2#3{\genspazio {#1}{#2}{#3}T0}

\def\L {\spazio L}
\def\H {\spazio H}
\def\W {\spazio W}

\def\C #1#2{C^{#1}([0,T];#2)}


\def\Lx #1{L^{#1}(\Omega)}
\def\Hx #1{H^{#1}(\Omega)}

\def\Luno{\Lx 1}
\def\Ldue{\Lx 2}

\def\Huno{\Hx 1}
\def\Hdue{\Hx 2}


\def\LQ #1{L^{#1}(Q)}


\let\theta\vartheta
\let\eps\varepsilon
\let\phi\varphi

\let\TeXchi\chi                         
\newbox\chibox
\setbox0 \hbox{\mathsurround0pt $\TeXchi$}
\setbox\chibox \hbox{\raise\dp0 \box 0 }
\def\chi{\copy\chibox}


\def\Betaeps{f_{1,\eps}}
\def\betaeps{f'_{1,\eps}}
\def\mueps{\mu_\eps}
\def\phieps{\phi_\eps}


\def\VD{V^*}
\def\phin{\phi_n}
\def\mun{\mu_n}

\def\phiz{\phi_0}

\def\muza{\mu_{0, \alpha}}
\def\nuza{\nu_{0, \alpha}}
\def\phia{\phi_{\alpha}}
\def\mua{\mu_{\alpha}}
\def\xia{\xi_{\alpha}}


\normalfont

\Begin{document}


\title{{\bf Hyperbolic relaxation of the\\ chemical potential in the \\viscous Cahn--Hilliard equation}}
\author{}
\date{}
\maketitle

\Bcenter
\vskip-1.3cm
{\large\bf Pierluigi Colli$^{(1)}$}\\
{\normalsize e-mail: {\tt pierluigi.colli@unipv.it}}\\[.4cm]
{\large\bf J\"urgen Sprekels$^{(2)}$}\\
{\normalsize e-mail: {\tt juergen.sprekels@wias-berlin.de}}\\[.6cm]

$^{(1)}$
{\small Dipartimento di Matematica ``F. Casorati'', Universit\`a di Pavia}\\
{\small and Research Associate at the IMATI -- C.N.R. Pavia}\\ 
{\small via Ferrata 5, 27100 Pavia, Italy}\\[.2cm]
$^{(2)}$
{\small Weierstrass Institute for Applied Analysis and Stochastics}\\
{\small Mohrenstra\ss e 39, 10117 Berlin, Germany}\\[.8cm]

\Ecenter

\begin{center}
\emph{In memory of Prof.~Dr.~Wolfgang Dreyer\\
with admiration, sympathy and friendship}
\end{center}

{
\Begin{abstract}\noindent
In this paper, we study a hyperbolic relaxation of the 
viscous Cahn--Hilliard system with zero Neumann boundary
conditions. In fact, we consider a relaxation term involving the 
second \juerg{time} derivative of the chemical potential in the first equation
of the system. We develop a well-posedness, continuous dependence
and regularity theory for the \juerg{initial-boundary} value problem. 
Moreover, we investigate the asymptotic behavior of the system as
the relaxation parameter tends to $0$ and prove the convergence
to the viscous Cahn--Hilliard system. 
\\[2mm]
{\bf Key words:}
Cahn--Hilliard system, hyperbolic relaxation, partial differential equations, initial-boundary value problem, well-posedness, continuous dependence, regularity, asymptotic convergence.
\normalfont
\\[2mm]
\noindent {\bf AMS (MOS) Subject Classification:}  
       35M33 
       35M87 
       35B40 
       37D35 
\End{abstract}
}
\salta

\pagestyle{myheadings}
\newcommand\testopari{\sc Colli \ --- \ Sprekels }
\newcommand\testodispari{\sc Hyperbolic relaxation of viscous Cahn--Hilliard equations}
\markboth{\testopari}{\testodispari} 


\finqui


\section{Introduction}
\label{Intro}
\setcounter{equation}{0}

In this paper, we deal with an initial-boundary value problem for a system of partial differential equations of viscous Cahn--Hilliard type, which in particular includes a hyperbolic relaxation term in the first equation.

The system is stated as follows:
\begin{align}
\label{ss1}
&\alpha \dtt\mu +  \dt\phi - \Delta \mu = 0 \quad \mbox{in }\,Q:= \Omega \times (0,T), \\
\label{ss2}
&\tau \dt \phi -\Delta \phi + f'(\phi) = \mu + g  \quad \mbox{in }\,Q,\\
\label{ss3}
&\dn \mu = \dn\phi = 0 \quad \mbox{on }\,\Sigma:= \partial \Omega \times (0,T),\\
\label{ss4}
&\mu(0)=\mu_0,\ \,  (\dt \mu) (0)=\nu_0 ,\ \, \phi(0)= \phi_0 \quad \mbox{in }\,\Omega,
\end{align}
\Accorpa\State ss1 ss4
where $\Omega\subset \erre^N$, $N\in\{1,2,3\}$, is a bounded and connected domain with smooth boundary  $ \partial\Omega$ and $T$ denotes some final time. We denote by ${\bf n}$ the unit outward 
normal to $ \partial\Omega$, with the associated outward normal derivative $\,\dn\,$. Note that $\dn$
appears in the \juerg{homogeneous Neumann} boundary conditions stated in \eqref{ss3} for both \juerg{the} variables 
$\mu$ and $\phi$. 

The equations \eqref{ss1}--\eqref{ss2} constitute a variation of the Cahn--Hilliard system (introduced in \cite{CH} and \juerg{first}
approached mathematically in \cite{EZ})
\begin{align}
\label{ch1}
&\dt\phi - \Delta \mu = 0 \quad \mbox{in }\,Q, \\
\label{ch2}
& -\Delta \phi + f'(\phi) = \mu + g  \quad \mbox{in }\,Q,
\end{align}
in which a viscosity term $\tau \dt \phi$ has been included in the second equation and where especially
the hyperbolic relaxation term $\alpha \dtt\mu$  has been added in the first equation. The viscous Cahn--Hilliard system
\begin{align}
\label{vch1}
&\dt\phi - \Delta \mu = 0 \quad \mbox{in }\,Q, \\
\label{vch2}
&\tau \dt\phi - \Delta \phi + f'(\phi) = \mu + g  \quad \mbox{in }\,Q,
\end{align}
is well known and \juerg{was} already investigated in a number of papers (see 
\cite{CM,CGM,CGS-Annali,CGS-JEPE, CGSCC,CoSpTr, CoGaMi, GiRoSi,GSS1,GS,KO} to quote some recent contributions), 
while to our knowledge an inertial term like 
$\alpha \dtt\mu$ in \eqref{ss1} is not common\revis{, but certainly} deserves to be examined. 
In this class of problems, the unknown 
functions $\,\phi\,$ and $\,\mu\,$ ususally stand for the \emph{order parameter}, which can represent a scaled
density of one of the involved phases, and the \emph{chemical potential}
associated with the phase separation process, respectively.

Moreover, $\,f' \,$ denotes the derivative 
(if it exists) of a  
double-well potential $f$, which in general  is split  into a (possibly nondifferentiable) convex part~$f_1$ 
 and a smooth and concave perturbation~$f_2$. Typical and physically significant examples for $\,f\,$ 
are the so-called {\em classical regular potential}, the {\em logarithmic double-well potential\/},
and the {\em double obstacle potential\/}, which are given, in this order,~by
\begin{align}
  & f_{\rm reg}(r) := \frac 14 \, (r^2-1)^2 \,,
  \quad r \in \erre, 
  \label{regpot}
  \\
  & f_{\rm log}(r):=\left\{\begin{array}{ll}
(1+r)\,\ln(1+r)+(1-r)\,\ln(1-r)-c_1r^2 &\quad\mbox{if }\,r\in(-1,1)\\
2\ln(2)-c_1 &\quad\mbox{if }\,r\in\{-1,1\}\\
+\infty&\quad\mbox{if }\,r\not\in [-1,1]
\end{array}, 
\right. 
  \label{logpot}
  \\[1mm]
  & f_{\rm 2obs}(r):=\left\{\begin{array}{ll}
c_2(1-r^2)  &\quad\mbox{if }\,r\in[-1,1]\\
+\infty&\quad\mbox{if }\,r\not\in [-1,1]
\end{array}. \right.
  \label{obspot}
\end{align}
Here, the constants $c_i$ in \eqref{logpot} and \eqref{obspot} satisfy
$c_1>1$ and $c_2>0$, so that $f_{\rm log}$ and $f_{\rm 2obs}$ are nonconvex.
Notice that for $f=f_{\rm log}$ the term $f'(\phi)$
occurring in \eqref{ss2} becomes singular as $\phi\searrow -1$ and $\phi\nearrow1$, which forces the order parameter $\phi$ to attain its values in the physically meaningful range~$(-1,1)$.  
In the nonsmooth case \eqref{obspot}, the convex part~$f_1$ is given by
the~indicator function of $[-1,1]$.
Accordingly, in such cases one has to replace the derivative of the convex part
by the subdifferential $\partial f_1$ and, consequently, to interpret \eqref{ss2} as a differential inclusion
or a variational inequality. 
We also note that $\tau$ is a \juerg{fixed} positive parameter (the viscosity coefficient), \juerg{while for the positive parameter
$\alpha$ we} will also discuss
the asymptotic convergence to $0$. We point out that in  \eqref{ss2} a known forcing term $g$ is present \juerg{that} may be 
interpreted as a direct or secondary control term which acts on the system. \juerg{In this connection}, we mention that optimal control problems
for viscous Cahn--Hilliard systems with a distributed control term involving $g$ \juerg{have recently been} treated in the paper \cite{CoSpTr}.

Some hyperbolic relaxations of the viscous Cahn--Hilliard system have been already considered and studied: let us mention the recent 
contributions~\cite{Bon,Chen,CMY,DMP,SZ,SK}. However, the available investigations are concerned with systems where the inertial 
term involves the phase variable $\phi$. In our case, the system \State\ couples a wave-type equation for $\mu$ \juerg{combined} with a source term given by $-\dt \phi$, with a semilinear parabolic equation in which the source term includes $\mu$. 

From the energetic viewpoint,
 there is a change with respect to the viscous (and nonviscous) Cahn--Hilliard equation. \juerg{To see this, let us for simplicity} argue now on the case $g=0$.  \juerg{Indeed}, for \eqref{vch1}--\eqref{vch2}, as well as for
 \eqref{ch1}--\eqref{ch2}, the basic energy estimate is obtained by testing \eqref{vch1} by $\mu$, \eqref{vch2} by $\dt \phi $, and \juerg{then adding and thus producing a cancellation of the terms containing the product $\mu \dt \phi $. Therefore}, one has 
 $$
 \int_0^t \!\!\int_\Omega |\nabla \mu|^2 + \tau  \int_0^t \!\!\int_\Omega |\dt\phi |^2  + \iO \Bigl( \frac12|\nabla \phi (t)|^2  
 +f(\phi(t))\Bigr) = \hbox{constant} $$
for $t\in [0,T]$, \juerg{where the first two terms are dissipative and the energy term is} given by the third one. The second term is missing in the case of the 
Cahn--Hilliard system~\eqref{ch1}--\eqref{ch2}, but the energy is the same and there is only one term 
for dissipation. On the other hand,
it turns out to be more difficult and involved to recover an energy estimate for \eqref{ss1}--\eqref{ss2}: as you will check in the sequel, our main  
estimate is constructed by testing \eqref{ss1} by $\dt \mu$ and the time derivative of \eqref{ss2} by $\dt\phi$, in order to have \juerg{a} cancellation of the terms containing the product $\dt \phi\, \dt \mu $. By integration, we then obtain 
\begin{align}
&\iO \Bigl( \frac \alpha 2|\dt \mu(t)|^2 + \frac 1 2 |\nabla \mu (t) |^2 +
\frac \tau 2|\dt \phi(t)|^2 \Bigr)\non\\
&\quad{} + \int_0^t \!\!\int_\Omega |\nabla(\dt\phi)|^2 +
 \int_0^t \!\!\int_\Omega f_{1}'' (\phin)|\dt\phi|^2 \non\\
&{} = \hbox{constant}\, {} -  \int_0^t \!\!\int_\Omega f_{2}'' (\phi)|\dt\phi|^2 \quad \hbox{for }\, t\in [0,T], 
\non
\end{align}
where the energy is now located in the first integral \juerg{in which neither the nonlinearity $f$ nor any of its derivatives occur}.  Moreover, note that 
the last term on the \lhs\ induces dissipation (as $ f_{1}'' $ is nonnegative), but on the \rhs\ the complementary term
may be positive and grow with respect to $t$, since $f_2$ is concave, in general. The viscous contribution in 
\eqref{ss2} is important here to control this term on the \rhs,  since the addendum $ \juerg{\frac {\tau} 2}|
\dt \phi(t)|^2$ is part in the energy. However, by our estimate we can proceed in the analysis and not only 
construct a well-posedness theory but also investigate the asymptotic behavior as the parameter $\alpha$ 
converges to $0$. 

This paper is dedicated to the memory of Wolfgang Dreyer, who recently passed away. The authors of this paper were fortunate to \juerg{have benefited} from Wolfgang's
friendship, as well as from his exceptional expertise and insight in Thermodynamics and Applied Mathematics. He was a brilliant and generous 
scientist who truly enjoyed engaging in scientific discussions with friends and colleagues. We both feel enriched by having known him and are grateful for the opportunity to have collaborated with him. Together, along with other colleagues, we co-authored the paper~\cite{BCDGSS} that explored the effects of phase separation driven by mechanical actions in tin/lead alloys, and where the corresponding system of partial differential equations already included equations of Cahn--Hilliard type. 

The paper is organized as follows. In the following section, we formulate the general assumptions and state the main results concerning 
the system~\State.  In Section~\ref{Exi}, we then prove the existence of \juerg{a} solution by using a double approximation based on a 
Yosida regularization of $\partial f_1$ and on a Faedo--Galerkin scheme. \juerg{This} proof requires the main analytical effort of this paper,
since it involves a number of estimates and two \juerg{passage-to-the-limit processes}. In Section~\ref{Cont}, we show the results on 
continuous dependence with respect to data and on the regularity \juerg{of} the solution: actually, three theorems are proved there. 
The final Section~\ref{Convergence} then brings the asymptotic results of the convergence of the system 
to the viscous Cahn--Hilliard system as $\alpha$ tends to $0$ and an estimate of the difference of solutions 
in terms of a precise rate of convergence.

We fix some notation.
For any Banach space $X$, we let $X^*$  denote 
its dual space, and $\norma{ \cpto }_X$ stands for the norm 
in $X$ and any power of $X$.   
For two Banach spaces $X$ and $Y$ that are both continuously embedded in some topological vector space~$Z$, the linear space
$X\cap Y$ is the Banach space equipped with its natural norm $\norma v_{X\cap Y}:=\norma v_X+\norma v_Y\,$  for $v\in X\cap Y$.
The standard Lebesgue and Sobolev spaces $L^p(\Omega)$ and $W^{m,p}(\Omega)$ are defined on $\Omega$ for
$1\le p\le\infty$ and $\,m \in \enne \cup \{0\}$. 
For the sake of 
convenience, we denote the norm of $\,L^p(\Omega)\,$ by $\,\|\,\cdot\,\|_p\,$ for $1\le p\le\infty$. 
If $p=2$, we employ the usual notation $H^m(\Omega):= W^{m,2}(\Omega)$. We also set
\begin{align*}
  & H := \Ldue , \quad V := \Huno, \quad W:=\bigl\{v\in \Hdue : \ \dn v =0 \, \hbox{ on }\, \Gamma \bigr\}.
  \end{align*}
Moreover, $\VD$~is the dual space of $V$,
and $\<\cpto,\cpto>$ stands for the duality pairing between $\VD$ and~$V$. We denote by $(\,\cdot\,,\,\cdot\,)$ the natural inner product in $\,H$.
As usual, $H$ is identified with a subspace of the dual space $\VD$ according to the identity
\begin{align*}
	\langle u,v\rangle =(u,v)
	\quad\mbox{for every $u\in H$ and $v\in V$}.
	\end{align*}
Note that $W \subset V \subset H\equiv H^* \subset \VD$ with dense and compact embeddings. 
About the constants used in the sequel for estimates, we adopt the rule that $\,C\,$ denotes any 
positive constant that depends only on the given data. The value of such generic constants $\,C\,$ 	may 
change from formula to formula or even within the lines of the same formula. Finally, the notation $C_\delta$ indicates a positive constant that additionally depends on the quantity $\delta$.   


\section{Main results}
\setcounter{equation}{0}
In this section, we formulate the general assumptions for the data of the system \State\ 
and state \juerg{existence}, continuous dependence, and regularity results.
First, \juerg{let us remark that the positive parameter $\alpha$ is not} listed in the assumptions below, 
since it is also involved \juerg{in} the related asymptotic analysis, and, consequently, we let 
$$ 0 < \alpha \leq 1 . $$ 
On the other hand, throughout the paper we suppose that
\begin{align}
& \tau>0  \, \hbox{ is a fixed constant.} 
\label{hyp1}
\end{align}
For the nonlinearity $f$ we generally assume that 
\begin{align}
& f=f_1+f_2, \quad\mbox{where}\non \\ 
&\quad{} f_1 :\erre\to [0,+\infty] \  \hbox{ is convex and lower semicontinuous with $f_1(0)=0$,}\non\\
&\quad{} f_2:\erre\to\erre \ \hbox{ has a Lipschitz continuous first derivative $f_2'$ on $\erre$.}\label{hyp2}
\end{align} 
A consequence of \eqref{hyp2} is that 
\Beq
\label{hyp3} 
\hbox{the subdifferential } \partial f_1 \hbox{ is maximal monotone 
in $\erre \times \erre$, with } 0 \in \partial f_1(0),
\Eeq
and an important requirement for the sequel is that 
\Beq
\label{hyp4} 
\hbox{the domain } D(\partial f_1) \hbox{  of $ \partial f_1 $ has a non-empty interior 
containing $0$.}
\Eeq
Note that these conditions are fulfilled in each of the cases considered in \eqref{regpot}, \eqref{logpot}, \eqref{obspot} with the domain $ D(\partial f_1)$ given by $\erre, \, (-1,1) , \, [-1,1]$, respectively.
From now onwards we use the symbol $\partial f_1^\circ (r)$ for 
the element of $\partial f_1 (r)$ (with $r\in D(\partial f_1)$) having minimum modulus, and 
we extend the notations $f_1$, $\partial f_1$, $D(\partial f_1)$, and $\partial f_1^\circ$
to the \juerg{corresponding} functionals and the operators induced on $L^2$ spaces. 

Also, \juerg{we assume} for the forcing term $g$ and the initial values $\mu_0 , \, \nu_0 , \, \phi_0$
that
\begin{align}
&g \in \H1H, \label{hyp5}\\
&\mu_0 \in V, \quad \nu_0 \in H, \label{hyp6}\\
& \phi_0 \in W \cap D(\partial f_1)  \ \hbox{ with } \ \partial f_1^\circ  (\phi_0) \in H. \label{hyp7}
\end{align} 
\Accorpa\Ipotesi hyp1 hyp8
Note that the 
condition $\phi_0 \in  W $ implies that $\phi_0\in C^0(\overline \Omega)$. Moreover, we require that
\begin{align}
&m_0 := \frac 1 {|\Omega|}\iO \phiz  \ \hbox{ lies in the interior of  } \ D(\partial f_1). \label{hyp8}
\end{align} 
\juerg{Here, $\,|\Omega|\,$ denotes the Lebesgue measure of $\,\Omega$,  and $m_0$ thus represents} the mean value of $\phiz$. 
In the following, we use the general notation $\overline  v$ to denote the mean value of a generic function $v\in\Luno$. If $v$ is in $V^*$,
then we can set
\Beq
  \overline v : = \frac 1 {|\Omega|} \, \langle v , 1 \rangle
  \label{defmean}
\Eeq
as well, noting that the constant function 1 is an element of $V$. Clearly, $\overline  v$ is the usual mean value of $v$ if $v\in H$. Note also that  $m_0 = \overline \phiz$. 

Let us now specify {our notion of} solution. We state the problem \State\ in a variational form.
We also introduce \juerg{an} additional variable $\xi$, which plays the role of $f_1'(\phi)$ in the case when the derivative of $f_1$ is replaced by a real subdifferential $\partial f_1$. 
Namely, the solution is a triple $(\mu,\phi,\xi)$ satisfying the regularity requirements
\begin{align}
  & \mu \in \W{2,\infty}{V^*}\cap \W{1,\infty}H \cap \L\infty V,
  \label{regmu}
  \\
  & \phi \in \W{1,\infty}H \cap \H1V \cap \L\infty W,
  \label{regphi}
  \\
  & \xi \in \L\infty H,
  \label{regxi}
\end{align}
\Accorpa\Regsoluz regmu regxi
and the following variational equations and initial conditions:
\begin{align}
  & \alpha \< \dtt\mu(t) , v > 
  + ( \dt\phi(t), v )
  + \iO \nabla\mu(t) \cdot \nabla v
  = 0
  \non
  \\
  & \quad \hbox{\aat\ and every $v\in V$},
  \label{prima}
  \\[1mm]
  &  \tau (\dt \phi(t) , v)
  + \iO \nabla\phi(t) \cdot \nabla v
  + ( \xi(t) +  f'_2 (\phi(t)), v) = (\mu(t) + g(t), v) 
  \non
  \\
  & \quad \hbox{\aat\ and every $v\in V$},
  \label{seconda}
  \\[2mm]
  &
  \xi \in \partial f_1 (\phi) \quad \aeQ ,  
  \label{terza}
  \\[2mm]
  &\mu(0)=\mu_0\, ,\quad  (\dt \mu) (0)=\nu_0 \, ,\quad  \phi(0)= \phi_0 \quad \hbox{ a.e.~in }{\Omega} \,.
  \label{cauchy}
\end{align}
\Accorpa\Pbl prima cauchy
Of course, in view of the regularities in \Regsoluz, the variational equality \eqref{seconda} is actually equivalent to an equation (cf.~\eqref{ss2}) plus the boundary condition for $\phi$ \juerg{which is already encoded in the fact that} 
$\phi \in \L\infty W$. Therefore, \eqref{seconda} can be replaced~by
\Beq 
  \tau \dt \phi -\Delta \phi + \xi + f_2'(\phi) = \mu + g 
  \quad \aeQ .
  \label{secondabvp}
\Eeq
On the contrary, the analogue equivalence for \eqref{prima} (cf.~\eqref{ss1})  would be true 
only if $\mu$ {were} more regular (cf. the first ad third term in~\eqref{prima}).

\Brem
\label{IntPbl}
Note that\juerg{, owing to the compactness of the embedding $W\subset C^0(\overline\Omega)$ for $N\le 3$, it follows from \cite[Sect.~8,~Cor.~4]{Simon}
and the regularity \eqref{regphi} that $\phi \in C^0(\overline Q).$  By the same token,
we have thanks to \eqref{regmu} that $\mu \in   \C1{V^*}\cap \C0H$, and, consequently,} $\dt \mu$ is at least weakly continuous from $[0,T]$ to $H$, which gives a meaning to the initial conditions in \eqref{cauchy}.  
\Erem

The next statement yields a well-posedness result for \Pbl .
\Bthm
\label{Wellposed}
Assume that \eqref{hyp1}--\eqref{hyp8} are fulfilled.
Then there exists a unique triple $(\mu,\phi,\xi)$,
with the regularity as in \Regsoluz, that solves problem \Pbl\
and satisfies the estimate
\begin{align}
&\alpha\norma{\mu}_{\W{2,\infty}{V^*}} + 
\alpha^{1/2}\norma{\mu}_{\W{1,\infty}H} 
+ \norma{\mu}_{\L\infty V} \non\\
&+ \norma{\phi}_{\W{1,\infty}H \cap \H1V \cap \L\infty W}
+ \norma\xi_{\L\infty H}
  \leq K_1
  \label{stimaK1}
\end{align}
for some constant $K_1>0$ that depends only on 
$\Omega$, $T$ and the data in \Ipotesi, but is independent of $\alpha$. 
\Ethm

The uniqueness property stated above is a consequence of the following continuous dependence result. 
Here, we use the notation
$$ 1*v (t) = \int_0^t v(s) ds ,  \quad \hbox{for } \, v\in L^1(0,T;\VD) \, \hbox{ at least}. $$ 

\Bthm
\label{Contdep}
Under the assumptions \eqref{hyp1}--\eqref{hyp4}, let $\, g_i, \mu_{0,i} , \, \nu_{0,i} , \, \phi_{0,i}$, $\, i=1,2$, be two sets of data satisfying \eqref{hyp5}--\eqref{hyp8}, and let $(\mu_i,\phi_i,\xi_i)$, $i=1,2$, denote any corresponding solutions \juerg{to} problem \Pbl\ with the
 regularity as in \Regsoluz. Then, the estimate  
\begin{align}
&\alpha^{1/2} \norma{\mu_1-\mu_2}_{\L\infty H}
+ \norma{\nabla (1* (\mu_1-\mu_2)) }_{\L\infty H}
+ \norma{\phi_1-\phi_2}_{\L\infty H\cap\L2V}
\non \\
&{}\leq K_2 \, \Bigl( \norma{g_1-g_2}_{\L2H} + \alpha^{1/2} \norma{\mu_{0,1} - \mu_{0,2}}_{H} +
\alpha^{1/2} \norma{\nu_{0,1} - \nu_{0,2}}_{H}     \Bigr)
\non \\
&\quad {}
+ K_2 \bigl(1+ \alpha^{- 1/2}\bigr) \norma{\phi_{0,1} - \phi_{0,2}}_{H}
  \label{stimaK2}
\end{align} 
holds true with a constant $K_2>0$ that depends only on $\Omega$, $T$, $\tau$, some Lipschitz constant for $f'_2$, and is independent of $\alpha$.
\Ethm

We observe that, by taking the same data in Theorem~\ref{Contdep}, the estimate \eqref{stimaK2} ensures uniqueness for the solution components $\mu $ and $\phi$ in the statement of Theorem~\ref{Wellposed}, while the uniqueness of $\xi$ results from \eqref{secondabvp} since the other terms in the equality are uniquely determined.  

On the basis of the estimates \eqref{stimaK1} and \eqref{stimaK2}, we are interested to investigate the asymptotic behavior of the problem~\Pbl\ as $\alpha \searrow 0$. This analysis will be developed in Section~\ref{Asy}, whereas now we discuss some further results that mostly depend on the values of $\alpha >0.$ A regularity result is stated below, under the additional assumption that 
\begin{align}
&\mu_0 \in W, \quad \nu_0 \in V.  \label{hyp9}
\end{align} 
\Bthm
\label{Regularity}
Assume that \eqref{hyp1}--\eqref{hyp8} and \eqref{hyp9} are fulfilled.
Then the unique solution  $(\mu,\phi,\xi)$ \juerg{to} problem \Pbl\ \juerg{satisfies} 
\begin{align}
  & \mu \in \W{2,\infty}{H}\cap \W{1,\infty}V \cap \L\infty W\,,
  \label{regmupiu}
\end{align}
and there exists a constant $K_3$, independent of $\alpha$, such \juerg{that}
\begin{align}
&\alpha\norma{\mu}_{\W{2,\infty}{H}} + 
\alpha^{1/2}\norma{\mu}_{\W{1,\infty}V} 
+ \norma{\mu}_{\L\infty W} 
  \leq K_3 \bigl(1+ \alpha^{- 1/2}\bigr). 
  \label{stimaK3}
\end{align}
\Ethm
Notice that, due to the \juerg{compactness of the embedding $W\subset C^0 (\overline{\Omega}) $, it follows from \cite[Sect.~8, Cor.~4]{Simon}} that 
$\mu \in \juerg{C^0(\overline Q)}$. 

For the next regularity and continuous dependence result we have to assume \juerg{further regularity for} $g$, that is,
\begin{align}
&g \in \LQ\infty, \label{hyp10}
\end{align}
and for the nonlinearity $f$. Namely, we suppose that the effective domain of $\partial f_1 $ is an open 
interval and that \juerg{the restriction of $f_1$ to} this interval is a smooth function. More precisely, we \juerg{assume that}  
\begin{align}
&D(\partial f_1 ) = (r_-,r_+), \quad  \hbox{with } \,-\infty\le r_-<0<r_+\le +\infty, \non \\
&\hbox{and the restriction of $f_1$ to $\,(r_-,r_+)\,$ belongs to $\,{C^2}(r_-,r_+)$}.
\label{hyp11}
\end {align}
\juerg{Then, for $r\in (r_-,r_+)$, the subdifferential $\partial f_1(r) $ reduces to the singleton
$\{f_1'(r)\} $}, and we require that
\begin{align}
	\text{$\lim_{r\searrow r_-} f_1'(r)=-\infty \, , \quad \lim_{r\nearrow r_+}f_1'(r)=+\infty$.}
	\label{hyp12}
\end{align}
Please note that both the potentials $f_{\rm reg}$ and $f_{\rm log}$ in 
  \eqref{regpot} and 
  \eqref{logpot} fulfill \eqref{hyp11}--\eqref{hyp12} with $\,(r_-,r_+)= \erre \,$ and $\,(r_-,r_+)= (-1,1),\,$ respectively. 

The so-called separation property and a refined continuous dependence result are stated as follows. 

\Bthm
\label{Separation}
Assume that \eqref{hyp1}--\eqref{hyp8} and \eqref{hyp9}--\eqref{hyp12} are fulfilled.
There exists two real numbers $r_*$ and~$r^*$, depending on $\alpha$ and on the structure of the system, such that
\begin{align}
		r_- < r_* \leq \phi(x,t)\leq r^* < r_+ \quad \text{ for \juerg{every} $ \,\,(x,t) \in \juerg{\overline Q}$}.
		\label{stima-sep}
\end{align}
Moreover, if for $i=1,2$ we let $(g_i, \mu_{0,i} , \, \nu_{0,i} , \, \phi_{0,i})$  be a set of data and 
$(\mu_i,\phi_i,\xi_i)$, with $ \xi_i = f_1'(\phi_i)$, denote the  corresponding solution \juerg{to} problem \Pbl, the estimate  
\begin{align}
&\norma{\mu_1-\mu_2}_{\H2\VD\cap\W{1,\infty} H \cap \L\infty V}
+ \norma{\phi_1-\phi_2}_{\H1H\cap\L\infty V\L2W}
\non \\
&{}\leq K_4 \, \Bigl( \norma{g_1-g_2}_{\L2H} + \norma{\mu_{0,1} - \mu_{0,2}}_{V} +
\norma{\nu_{0,1} - \nu_{0,2}}_{H} + \norma{\phi_{0,1} - \phi_{0,2}}_{V}   \Bigr)
  \label{stimaK4}
\end{align} 
holds true for some constant $K_4>0$ that depends on $\alpha$ and on the structure of the system. 
\Ethm

\Brem
Please note that in the case when the domain $D(\partial f_1 )$ is the entire real line, i.e., if
$r_- = -\infty$ and $ r_+ = +\infty$, then the property~\eqref{stima-sep} is direcly ensured by 
the estimate~\eqref{stimaK1}, since \juerg{$\phi$ is bounded in $C^0(\overline Q)$};
therefore, if  $D(\partial f_1 )= \erre$, then the additional regularity assumptions \eqref{hyp9} 
and \eqref{hyp10} are not needed to prove~\eqref{stima-sep}.
\Erem

\Brem
It would be interesting to investigate the system \eqref{ss1}--\eqref{ss4} from the viewpoint of an optimal control problem, with the distributed control located in the source term $g$ in equation \eqref{ss2}. Thus,
in order to discuss differentiability properties and optimality conditions, it could be important to deal with smoother data and a 
smooth nonlinearity $f$, and with the control $g$ lying in a control box in $\LQ\infty$ (cf.~\eqref{hyp10}). 
In this framework, stronger stability and continuous depencence estimates 
\juerg{can possibly be} derived for the system. \juerg{However, the estimates}
 \eqref{stimaK1}, \eqref{stimaK3} and \eqref{stimaK4} are already a good starting point in that direction. 
\Erem

\section{Existence of solutions}
\label{Exi}
\setcounter{equation}{0}

In this section, we are going to prove the existence result for the problem \Pbl, by constructing a solution 
$\,(\mu, \phi,\xi)\,$ that satisfies \Regsoluz. We adopt two levels of approximation: \juerg{at first,} we replace the 
subdifferential $\partial f_1$ in \eqref{terza} by the derivative of the Moreau--Yosida regularization $\Betaeps$ 
of $f_1$, depending on a parameter  $\varepsilon\in (0,1)$; then, \juerg{we apply a Faedo--Galerkin scheme to the resulting approximate system}. 

To begin with, we consider for every $\varepsilon\in (0,1)$ the Moreau--Yosida regularization $\Betaeps$ of $f_1$, that is (see, e.g., \cite{Barbu, Brezis}),
\begin{align*}
	& \Betaeps(r)
	:=\inf_{s \in \mathbb{R}}\left\{ \frac{1}{2\varepsilon } |r-s|^2
	+f_1(s) \right\} 
	= 
	\frac{1}{2\varepsilon } 
	\bigl| r-J_\varepsilon  (r) \bigr|^2+f_1 (J_\varepsilon (r) )
	= \int_{0}^{r} \betaeps (s)ds,
\end{align*}
where $\betaeps :\mathbb{R} \to \mathbb{R}$ \juerg{and} the associated resolvent operator
$J_\varepsilon $ \juerg{are} given 
by 
\begin{align*}
	& \betaeps (r)
	:= \frac{1}{\varepsilon } ( r-J_\varepsilon (r) ), 
	\quad 
	J_\varepsilon (r) 
	:=(I+\varepsilon  \partial f_1 )^{-1} (r)\juerg{, \quad\mbox{for all \,$r\in\erre$,}}
\end{align*}
\juerg{with} $I$ denoting the identity operator. Note that the derivative $\betaeps$ turns out to be a regularization of 
the graph $\partial f_1$. Indeed, $\betaeps$ and
$\Betaeps$ fullill, for all $0<\eps< 1$ (see, e.g., \cite[pp.~28 and~39]{Brezis}), 
\begin{align}
  & \hbox{$\betaeps :\erre \to \erre $ is monotone and \Lip\ continuous} 
	\nonumber\\
	&\quad \hbox{with Lipschitz constant } 1/\eps, \mbox{ and it holds} \,\betaeps(0)=0,
  \label{funoeps1}
  \\
  & |\betaeps(r)| \leq |\partial f_1^\circ (r)|
  \quad \hbox{for every $r\in D(\partial f_1)$},
  \label{funoeps2}
  \\
  & 0 \leq \Betaeps(r)  \leq f_1(r)
  \quad \hbox{for every $r\in\erre$}.
  \label{funoeps3}
\end{align}
\Accorpa\Propbetaeps funoeps1 funoeps3

As for the second approximation, we employ a Faedo--Galerkin discrete scheme using a special basis. 
To this end, we take the eigenvalues $\,\{\lambda_j\}_{j\in\enne}\,$ of the eigenvalue problem 
$$-\Delta v=\lambda v \quad\mbox{in }\,\Omega, \qquad \dn v=0 \quad\mbox{on }\,\partial \Omega,$$
and let $\,\{e_j\}_{j\in\enne}\subset W\,$ be the associated eigenfunctions, normalized by $\,\|e_j\|_H=1$, $j\in\enne$. Then, we have that
\begin{align*}
&0=\lambda_1<\lambda_2\le \ldots, \qquad \lim_{j\to\infty}\lambda_j=+\infty,\\
&\iO e_je_k=\iO \nabla e_j\cdot\nabla e_k=0\quad \mbox{for $\,j\not= k$},
\end{align*}
and we note that $e_1$ is just the constant function $|\Omega|^{-1/2}.$
We then define the $n$-dimensional spaces $\,V_n:={\rm span}\{e_1,\ldots,e_n\}$ for $\,n\in\enne$, where $V_1$ \juerg{is just} 
the space of constant functions on $\,\Omega$. 
It is well known that the union of these spaces is dense in both $\,H\,$ 
and~$\,V$. 

The approximating $n$-dimensional problem is stated as follows:
find functions
\begin{equation}
\label{discrete}
\mu_n(x,t)=\sum_{j=1}^n \mu_{nj}(t)e_j(x),\quad \phi_n(x,t)=\sum_{j=1}^n \phi_{nj}(t)e_j(x),
\end{equation}
such that 
\begin{align}
  & \alpha ( \dtt\mu_n(t) , v ) 
  + ( \dt\phi_n(t), v) 
  + \iO \nabla\mu_n(t) \cdot \nabla v
  = 0
  \non
  \\
  & \quad \hbox{for all $t\in [0,T]$ and every $v\in V_n$},
  \label{ss1n}
  \\[1mm]
  &  \tau (\dt \phi_n (t), v)
  + \iO \nabla\phi_n (t) \cdot \nabla v
  + ( \betaeps (\phi_n(t))+  f'_2 (\phi_n(t)), v) = 
  (\mu_n(t) + g(t), v) 
  \non
  \\
  & \quad \hbox{for all $t\in [0,T]$ and every $v\in V_n$},
  \label{ss2n}
  \\[2mm]
  &\mu_n(0)=P_n(\mu_0)\, ,\quad  (\dt \mu_n) (0)=P_n(\nu_0) \, ,\quad  \phi_n(0)= P_n(\phi_0) \quad \hbox{ a.e.~in }{\Omega} \,,
  \label{ss4n}
\end{align}
where $P_n$ denotes the $H$-orthogonal projection onto $V_n$. Then $P_n(v)=\sum_{j=1}^n (v,e_j)e_j$ for every $v\in H$, and we have  
(see, e.g., \cite[formula~(3.14)]{CGSS6})
\begin{equation}
\label{Fiete}
\|P_n (v)\|_Y\,\le\,C_\Omega \|v\|_Y \quad\mbox{for every $\,v\in Y$, where }\,Y\in\{H,V,W\},
\end{equation}
for some constant $C_\Omega>0$ depending only on $\,\Omega$. By comparing \eqref{ss1n}--\eqref{ss4n} with \Pbl, note that the inclusion \eqref{terza} present in \Pbl\ is not reproduced in
\eqref{ss1n}--\eqref{ss4n}, since the role of the $\xi$ variable is now played by $\betaeps (\phi_n)$, written as it is, in  \eqref{ss2n}.

Next, we take $v=e_k$ in all of the equations \eqref{ss1n}--\eqref{ss4n}, for $k=1,\ldots,n$, obtaining the system
\begin{align}
\label{ODE1}
&\alpha \frac{d^2}{dt^2} \mu_{nk} + \frac{d}{dt} 
\phi_{nk}+\lambda_k\,\mu_{nk}=0 \quad \mbox{in }\,(0,T),\\
\label{ODE2}
&\tau\frac{d}{dt} \phi_{nk}+\lambda_k\,\phi_{nk}
+(\betaeps(\phin) + f_2 (\phin),e_k)= \mu_{nk} +
(g, e_k) 
\quad \mbox{in }\,(0,T),\\
\label{ODE4}
& \mu_{nk}(0)=(\mu_0,e_k), \quad \frac{d}{dt}\mu_{nk}(0)=(\nu_0,e_k), \quad \phi_{nk}(0)=(\phiz,e_k).&&{}
\end{align}
Then we have to deal with a Cauchy problem for a system of ordinary differential equations\juerg{, which is of second order}
in the variables $\mu_{nk}$ 
and \juerg{of} first order in the variables $\phi_{nk}$. This system is set in explicit form and offers Lipschitz continuous nonlinearities 
and source terms $(g, e_k) $ in $H^1(0,T)$ (cf.~\eqref{hyp5}). 
By Carath\'eodory's theorem, the Cauchy problem \eqref{ODE1}--\eqref{ODE4} has a unique solution expressed by $\mu_{nk}, \, \phi_{nk}$, with $ \mu_{nk} \in H^3(0,T)$ and $ \phi_{nk} \in H^2(0,T)$, for $k=1,\ldots,n$. On account of \eqref{ODE1}--\eqref{ODE2} and \eqref{discrete},  this solution 
uniquely determines a pair $(\mun,\phin)\in H^3(0,T;V_n)\times H^2(0,T;V_n)$ that solves \eqref{ss1n}--\eqref{ss4n}. 

We now derive a series of a priori
estimates for the finite-dimensional approximations. 
In the following, $C>0$ denotes constants that may depend on the data of 
the state system, but \juerg{are} independent of $n\in\enne,\, \eps \in (0,1)$ and $\alpha \in (0,1]$.

\step
First estimate

The aim is taking the time derivative of \eqref{ss2n} and then 
testing by $\dt \phi_n$. In addition, we add the resulting 
equality to \eqref{ss1n}
where we choose $v= \dt\mu_n $. By this approach we obtain the 
cancellation of two terms. Next, we integrate with respect to time 
and, in order to recover our estimate, we have to control the $H$-norm of the initial value $\dt \phi_n (0) $. Taking $t=0 $ and 
$v=\dt\phi_n(0)$ in \eqref{ss2n}, by \eqref{ss4n} we easily infer that 
\begin{align*}
&\tau \norma{\dt \phi_n (0)}^2_H = (\Delta \phi_n (0)
- \betaeps (\phi_n(0)) - f'_2 (\phi_n(0)) + \mu_n(0) + g(0), \dt \phi_n (0))\\
&\quad{}\leq \frac\tau 2 \norma{\dt \phi_n (0)}^2_H + \frac1{2\tau} 
\norma{\Delta P_n(\phi_0)
- \betaeps (P_n(\phi_0)) - f'_2 (P_n(\phi_0)) + P_n(\mu_0) + g(0)}^2_H ,
\end{align*}
thanks to the Schwarz and Young inequalities. By virtue of the assumptions \eqref{hyp5}--\eqref{hyp7} on $g$ and the initial data,
and of the property~\eqref{Fiete}, it turns out that there is a constant $C_\eps$, depending only on the data and on $\eps$, such that  
\begin{align}
\label{esti1} 
&\tau \norma{\dt \phi_n (0)}_H  
\leq \norma{\Delta P_n(\phi_0) - \betaeps (P_n(\phi_0)) - f'_2 (P_n(\phi_0)) + P_n(\mu_0) + g(0)}_H\le C_\eps\,.
\end{align}
The dependence on $\eps$ follows from the Lipschitz continuity of $\betaeps$ with constant $1/\eps$ (cf.~\eqref{funoeps1}), while $f'_2$ is Lipschitz continuous independently of $\eps$ 
(see~\eqref{hyp2}). 

Now, we can perform the computation described above and deduce that
\begin{align}
&\frac \alpha 2\|\dt \mun(t)\|_H^2 + \frac 1 2\|\nabla \mun(t)\|_H^2 +
\frac \tau 2\|\dt \phin(t)\|_H^2 \non\\
&\quad{} +\iint_{Q_t} |\nabla(\dt\phin)|^2 +
\iint_{Q_t} f_{1,\eps}'' (\phin)|\dt\phin|^2 \non\\
& \leq\,
\frac \alpha 2\|P_n (\nu_0)\|_H^2 +
\frac 1 2\|\nabla P_n(\mu_0) \|_H^2 
\non\\
&\quad{}+ \frac1{2\tau} \norma{\Delta P_n(\phi_0) - \betaeps (P_n(\phi_0)) - f'_2 (P_n(\phi_0)) + P_n(\mu_0) + g(0)}^2_H \non\\
&\quad{} -\iint_{Q_t} f_2'' (\phin))|\dt\phin|^2 + \iint_{Q_t} \dt g \, \dt\phin,
\label{help1}
\end{align}
where we have used the notation
$$ Q_t:= \Omega \times (0,t), \quad t\in (0,T].$$
By the monotonicity of $\betaeps$, the last term on the \lhs\ of 
\eqref{help1} 
is nonnegative. Moreover, owing to \eqref{hyp2}, we have that  $\,|f_2''(\phin)||\dt\phin|^2 \le C |\dt\phin|^2 \,$  
a.e.~in~$Q_t$. Then, in view of \eqref{hyp5}, \eqref{hyp6}\juerg{, \eqref{Fiete},
Young's inequality, and Gronwall's lemma,} it is straightforward to infer 
from~\eqref{help1} that
\begin{align}
&\alpha^{1/2} \|\dt \mun\|_{\L\infty H} + \|\nabla \mun\|_{\L\infty H} +
\|\dt \phin\|_{\L\infty H} +\|\nabla(\dt\phin)\|_{\L2 H} \non\\
& \leq\,
C \big( \|\nu_0\|_H +
\|\mu_0 \|_V + \norma{ \dt g }_{\L2 H} 
\non \\
&\qquad\quad{}
+ \norma{\Delta P_n(\phi_0) - \betaeps (P_n(\phi_0)) - f'_2 (P_n(\phi_0)) + P_n(\mu_0) + g(0)}_H \, \big),
\label{pier1}
\end{align}
for some constant $C$ depending only on data, as $ 0 <\alpha \le 1 $. Therefore, recalling the initial conditions in \eqref{ss4n}, since 
$$
\mu_n(t) = P_n(\mu_0) + \int_0^t \dt \mu_n (s) ds , \quad
\phi_n(t) = P_n(\phi_0) + \int_0^t \dt \phi_n (s) ds  \quad
\hbox{ for all } t\in [0,T], 
$$
we easily conclude from \eqref{esti1} that 
\begin{equation}
\label{esti2}
\alpha^{1/2} \|\mun\|_{\W{1,\infty}H}
+  \|\nabla \mun \|_{\L{\infty}H}
+ \|\phin\|_{\W{1,\infty}H\cap \H1V}\,\le\,C_\eps\,.
\end{equation}

\step 
Complementary estimates

Taking now \revis{in \eqref{ss1n} and using the othogonality of the projection operator $P_n$ and} \eqref{Fiete}, we have that 
\begin{align}
&\alpha \< \dtt\mu_n (t) , v > =  
\alpha ( \dtt\mu_n (t) , P_n (v) )  + \alpha ( \dtt\mu_n (t) , v - P_n (v) )
\non\\
&\leq |\alpha ( \dtt\mu_n (t) , P_n (v) ) |
\leq \biggl|  ( \dt\phi_n (t), P_n(v) )  + \iO \nabla\mu_n (t) \cdot \nabla P_n(v)
\biggl| \non \\
&\leq C \bigl(\|\dt \phin\|_{\L{\infty}H} +  \|\nabla \mun \|_{\L{\infty}H} \bigr)  \norma{v}_V \quad \hbox{for a.e. } t\in (0,T), 
\label{estP1}
\end{align}
\juerg{so that} from \eqref{esti2} it clearly follows that
\begin{equation}
\label{estP2}
\alpha \|\dtt \mun\|_{\L{\infty}{\VD}} \,\le\,C_\eps\,.
\end{equation}
In addition, we can take $v= -\Delta (\phi_n(t) ) $ in \eqref{ss2n} and integrate by parts in some term. With the help of Young's inequality and the Lipschitz continuity of $f_2'$, we obtain
\begin{align}
& \juerg{\|\Delta \phin(t)\|^2_H} + 
\iO f_{1,\eps}''(\phi_n (t) )|\nabla \juerg{\phin(t)}|^2 \non\\
& = (\tau \dt \phi_n (t)+ f_2' (\phin (t)) - g(t),
\Delta \phin(t)) + \iO \nabla \mu_n (t) \cdot \nabla \phin(t) \non\\
& \leq \frac 1 2 \juerg{\|\Delta \phin(t)\|^2_H}  + 
C \bigl(1+ \|\phin\|_{\W{1,\infty}H}^2  +  \|g\|_{\L{\infty}H}^2 \bigr)
\non \\ 
&\quad{}+ \|\nabla \mun \|_{\L{\infty}H} \|\nabla \phin \|_{\L{\infty}H}
\quad \hbox{for a.e. } t\in (0,T),
\label{estP3}
\end{align}
where the second term in the first line is nonnegative due to \eqref{funoeps1}. Consequently, from \eqref{esti2} and the elliptic regularity theory,
we find that 
\begin{equation}
\label{estP4}
 \| \Delta \phin \|_{\L{\infty}H} +  \| \phin \|_{\L{\infty}W} \,\le\,C_\eps\,.
\end{equation}

\step
Passage to the limit in the Faedo--Galerkin scheme

By virtue of the uniform estimates shown above, there exists a pair $(\mueps,\phieps)$ such that (possibly on a subsequence, 
which is again labeled by $n\in\enne$)
\begin{align}
\label{conmu}
&\mun\to\mueps \quad \mbox{weakly star in }\,\W{2,\infty}{V^*}\cap \W{1,\infty}H \cap \L\infty V,\\
\label{conphi}
&\phin\to\phieps \quad \mbox{weakly star in } \,W^{1,\infty}(0,T;H)\cap H^1(0,T;V)\cap L^\infty(0,T;W).
\end{align}
By \eqref{conmu} and the compact embeddings $V\subset H \subset \VD$, it follows from \cite[Sect.~8,~Cor.~4]{Simon} that
\begin{equation}
\label{conmu1}
\mun\to\mueps\quad \mbox{strongly in }\,\C1{V^*}\cap \C0 H\,. 
\end{equation}
Moreover, owing to the compactness of the embeddings $W\subset C^0(\overline\Omega)$ 
and $W\subset V$, it turns out that 
\begin{equation}
\label{conny1}
\phin\to \phieps\quad\mbox{strongly in }\,C^0(\overline Q)\cap \C0V \,,
\end{equation}
whence, by the Lipschitz continuity of $\betaeps$ and $f'_2$, we deduce that 
\begin{equation}
\label{conny2}
\betaeps(\phin) + f_2'(\phin) \to \betaeps(\phieps) + f_2'(\phieps) \quad\mbox{strongly in }\,C^0(\overline Q)\cap \C0H\,,
\end{equation}
at least. Then, taking first $v\in V_k$ with $k\leq n$ in \eqref{ss1n}--\eqref{ss2n},
and passing to the limit as $n\to \infty $, it is not difficult to infer that 
\begin{align}
  & \alpha \< \dtt\mueps , v >
  + ( \dt\phieps, v) 
  + \iO \nabla\mueps \cdot \nabla v
  = 0 \quad \hbox{a.e. in $(0,T)$},
  \label{lim1}\\
  &  \tau (\dt \phieps, v)
  + \iO \nabla\phieps \cdot \nabla v
  + ( \betaeps (\phieps) +  f'_2 (\phieps), v) = 
  (\mueps + g, v) 
  \quad \hbox{a.e. in $(0,T)$},
  \label{lim2}
\end{align}
\juerg{at} first for all $v\in\cup_{k\in\enne} V_k$, \juerg{and} then, by density, for all $v\in V$. Note that, 
due to the regularity of $\phieps $, the second variational equality can be equivalently rewritten as 
\begin{equation}
\label{lim2-bis}
\tau \dt \phieps - \Delta \phieps +  \betaeps (\phieps) +  f'_2 (\phieps) = 
\mueps + g  \quad \hbox{a.e. in $Q$}, 
\end{equation}
\juerg{where it is} understood that $\phieps $ satisfies the boundary condition $\dn \phieps = 0$
a.e.~on $\Sigma$, on account of $\phieps \in\L\infty W$.

Thanks to \eqref{conmu1} and \eqref{conny1}, we can pass to the limit as $n\to\infty$ also in the initial conditions \eqref{ss4n} and find that
\begin {equation}
\mueps(0)=\mu_0\, ,\quad  (\dt \mueps) (0)= \nu_0  \, ,\quad  \phieps(0)= \phi_0\,, \quad \hbox{ a.e.~in }{\Omega} ,
  \label{iceps}
\end{equation}
since (see~\eqref{hyp6}--\eqref{hyp7})
$$ P_n(\mu_0) \to \mu_0, \quad P_n(\nu_0) \to \nu_0, \quad P_n(\phi_0) \to \phi_0, \quad \hbox{ strongly in } \, V, \, H, \, W ,\, \hbox{ respectively,} $$ 
and $\dt \mueps $ is weakly continuous from $[0,T] $  to $H$. Also, we can invoke the weak star lower semicontinuity of norms and pass to the limit 
in~\eqref{pier1} to derive the inequality 
\begin{align}
&\alpha^{1/2} \|\dt \mueps\|_{\L\infty H} + \|\nabla \mueps\|_{\L\infty H} +
\|\dt \phieps\|_{\L\infty H} +\|\nabla(\dt\phieps)\|_{\L2 H} \non\\
& \leq\,
C \big( \|\nu_0\|_H +
\|\mu_0 \|_V + \norma{ \dt g }_{\L2 H} 
\non \\
&\qquad\quad{}
+ \norma{\Delta \phi_0 - \betaeps (\phi_0) - f'_2 (\phi_0) + \mu_0 + g(0)}_H \, \big). 
\label{pier1-bis}
\end{align}
Hence, recalling \eqref{hyp7} and \eqref{funoeps2}, it turns out especially that $\norma{\betaeps (\phi_0)}_H $ is bounded independently of $\eps$, which \juerg{implies} that the complete \rhs\ of \eqref{pier1-bis} is uniformly bounded. Then, arguing as before,
we can improve \eqref{esti2} for $(\mueps, \phieps)$ and recover \juerg{that}  
\begin{equation}
\label{estP5}
\alpha^{1/2} \|\mueps\|_{\W{1,\infty}H}
+  \|\nabla \mueps \|_{\L{\infty}H}
+ \|\phieps\|_{\W{1,\infty}H\cap \H1V}\,\le\,C\,.
\end{equation}
As a consequence, by repeating the arguments in \eqref{estP1} and 
\eqref{estP3} for $\mueps $ and $\phieps$, we easily find out that 
\begin{equation}
\label{estP6}
\alpha \|\dtt \mueps\|_{\L{\infty}\VD}
+  \|\Delta  \phieps \|_{\L{\infty}H}
+ \|\phieps\|_{\L{\infty}W}\,\le\,C
\end{equation}
for some constant $C$ \juerg{which is} independent of \juerg{both} $\eps \in (0, 1)$ and $\alpha \in (0, 1]$.

\step 
Further estimates

We insert the constant function $v=1/|\Omega|$ in \eqref{lim1} and\revis{, in view of the definition~\eqref{defmean} of the mean value,} deduce that $\,
\overline{\dt (\alpha \dt\mueps + \phieps) }=0\,$ a.e. in
$(0,T)$. Hence, \revis{by virtue of \eqref{iceps}, \eqref{hyp8} and \eqref{defmean} again,}  it is straightforward to obtain
\begin{equation}
\label{pier2}
\alpha \,\overline{\dt\mueps}(t) + \overline{\phieps} (t) =
\alpha \, \overline{\nu_0} + m_0
\quad \hbox{for all $t\in [0,T]$.} 
\end{equation}
Now, we take $v= \phieps (t) - m_0 $ in \eqref{lim2} and, without integrating with respect to time, we have that 
\begin{align}
& \iO | \nabla (\phieps (t) - m_0 ) |^2
+ ( \betaeps (\phieps(t)), \phieps (t) - m_0 ) \non\\
&=  -(\tau \dt \phieps (t) +  f'_2 (\phieps(t)) -g(t), \phieps (t) - m_0) + (\mueps(t) , \phieps (t) - m_0 ). 
\label{pier3}
\end{align}
Now, in view of the properties \eqref{hyp2}--\eqref{hyp4} 
of $f_1$ and $\partial f_1$, and on account of 
\eqref{hyp8}, it turns out that there exist two positive constants $\delta_0$ and $C_0$, independent of $\eps$, such that 
\begin{equation}
\label{Zelik}
\betaeps(r)(r-m_0)\ge \delta_0|\betaeps(r)|-C_0 \quad\mbox{for every $\,r\in \erre$}.
\end{equation}%
For this property we refer to~\cite[Appendix, Prop.~A.1]{MiZe}
and also to the detailed proof given in~\cite[p.~908]{GiMiSchi}.
Applying \eqref{Zelik} to the second term in the \lhs\ of \eqref{pier3}, due to \eqref{hyp5} on $g$ and to the Lipschitz continuity of $f_2'$,
we infer that
\begin{align} 
\delta_0| \betaeps (\phieps(t))| 
\leq C \bigl(1+ \|\phieps\|_{\W{1,\infty}H}^2\bigr)
+ (\mueps(t) , \phieps (t) - m_0 ). 
\label{pier4}
\end{align}
As for the last term in \eqref{pier4}, thanks to \eqref{pier2} and the Poincar\'e--Wirtinger inequality, we can argue as follows:
\begin{align} 
&(\mueps(t) , \phieps (t) - m_0 ) \non\\
&=
(\mueps(t) , \alpha \,\dt\mueps(t) + \phieps (t) - \alpha \, \overline{\nu_0} - m_0 ) +
(\mueps(t) ,  \alpha \, \overline{\nu_0} - \alpha \,\dt\mueps(t)) 
\non\\
&= (\mueps(t) - \overline{\mueps} (t) , \alpha \,\dt\mueps(t) + \phieps (t) - \alpha \, \overline{\nu_0} - m_0 )
+ \alpha( \mueps(t) ,  \overline{\nu_0} - \,\dt\mueps(t))
\non\\
&\leq  C \|\nabla \mueps \|_{\L{\infty}H} \bigl(\alpha^{1/2} \|\dt\mueps\|_{\L{\infty}H} + \|\phieps \|_{\L{\infty}H} +1\bigr) 
\non\\
&\quad{} + C \bigl( \alpha \|\mueps\|_{\W{1,\infty}H}^2 +1 \bigr)
\label{pier5}
\end{align}
since $0<\alpha\leq 1$. 
Hence, by \eqref{pier4}, \eqref{pier5}, and \eqref{estP5}, we can conclude that 
\begin{equation}
\label{esti3}
\|\betaeps (\phieps)  \|_{\L{\infty}{L^1(\Omega)}}
\,\le\,C \, . 
\end{equation}
Next, we choose the constant function $\,v=1\,$ in \eqref{lim2}.  We obtain, for a.e. $t\in(0,T)$,
\begin{align}
\label{Uwe}
\iO \betaeps (\phieps(t)) + \iO (\tau \dt \phieps + f'_2(\phieps (t)) - g(t)) \,=\, |\Omega|\,\overline{\mueps}(t).
\end{align}
Owing to~\eqref{esti3} and \eqref{estP5}, both summands on the \lhs\ are bounded in $L^\infty(0,T)$. 
Then we infer that
\begin{equation}
\label{estP7}
\|\overline{\mueps}\|_{L^\infty(0,T)}\le C.
\end{equation}
Moreover, using again the Poincar\`e--Wirtinger inequality, we have that 
\begin{align}
\label{esti4}
&\|\mueps\|_{\L{\infty}H} \leq \|\mueps - \overline{\mueps} \|_{\L{\infty}H} +  C\|\overline{\mueps}\|_{L^\infty(0,T)} 
\non \\
&\quad{}
\le C \|\nabla \mueps \|_{\L{\infty}H} + C \leq C. 
\end{align}
Now we can go back to \eqref{lim2} or, better, \juerg{to} \eqref{lim2-bis} and compare the terms in the equation: 
from \eqref{estP5}, \eqref{estP6} and \eqref{hyp5} it follows that 
$\dt \phieps, \, \Delta \phieps, \,  f'_2 (\phieps), \,  \mueps, \, g $ are all uniformly bounded in $\L\infty H$, whence 
\begin{equation}
\label{esti5}
\|\betaeps (\phieps)  \|_{\L{\infty}H}
\,\le\,C \, . 
\end{equation}

\step
Passage to the limit as $\eps \to 0\, $

Thanks to the uniform estimates \eqref{estP5}, \eqref{estP6}, \eqref{esti4}, and \eqref{esti5},
 it follows that there is a triple $(\mu,\phi, \xi)$ such that,
for some subsequence $\eps_k$ tending to $0$, it holds  
\begin{align}
\label{conv1}
&\mu_{\eps_k} \to\mu \quad \mbox{weakly star in }\,\W{2,\infty}{V^*}\cap \W{1,\infty}H \cap \L\infty V,\\
\label{conv2}
&\phi_{\eps_k}\to\phi \quad \mbox{weakly star in } \,W^{1,\infty}(0,T;H)\cap H^1(0,T;V)\cap L^\infty(0,T;W),\\
\label{conv3}
&f'_{1,\eps_k} (\phi_{\eps_k}) \to\xi \quad \mbox{weakly star in } L^\infty(0,T;H).
\end{align}
As argued in the previous limit procedure~(cf.~\eqref{conmu1}--\eqref{conny2}), by compactness, in particular exploiting 
\cite[Sect.~8,~Cor.~4]{Simon}, and by the Lipschitz continuity of $f_2'$, we have that 
\begin{align}
\label{conv4}
&\mu_{\eps_k} \to\mu \quad\mbox{strongly in }\,\C1{V^*}\cap \C0 H\, ,
\\
\label{conv5}
&\phi_{\eps_k}\to\phi \quad \mbox{strongly in } \,C^0(\overline Q)\cap \C0V\, ,
\\
\label{conv6}
&f_2' (\phi_{\eps_k}) \to f_2'(\phi) \quad \mbox{strongly in } \,C^0(\overline Q)\cap \C0H\, .
\end{align}
Then we can pass to the limit as $\eps_k\to 0 $ in \eqref{lim1} and \eqref{lim2} by 
finding \eqref{prima} and \eqref{seconda}, respectively. Moreover, the initial conditions \eqref{cauchy} follow from \eqref{iceps}. It remains to check \eqref{terza}: but, since the extension of $\partial f_1 $ to $\L2H$ is a maximal monotone operator and $\betaeps$ denotes its Yosida approximation, 
\juerg{and since we have that} 
$$ \limsup_{k,n\to \infty} \int_0^T (f'_{1,\eps_k}(\phi_{\eps_k}(t)) - f'_{1,\eps_n}(\phi_{\eps_n}(t)), \phi_{\eps_k}(t)- \phi_{\eps_n}(t)) dt =0 $$
due to the weak convergence of $f'_{1,\eps_k} (\phi_{\eps_k})$ to $\xi$ in $\L2H$ and the strong convergence of $\phi_{\eps_k}$ to $\phi$ in $\L2H$, we can apply \cite[Prop.~2.2, p.~38]{Barbu} and recover the inclusion $\xi \in \partial f_1 (\phi)$ 
in $\L2H$ and almost everywhere in $Q$. 
 
In conclusion, we note that the triple $(\mu,\phi, \xi)$ found by the limit procedure is actually the unique solution of the problem \eqref{prima}--\eqref{cauchy}, on account of the 
continuous dependence result,  
and that the estimate \eqref{stimaK1} follows easily from 
the uniform bounds in \eqref{estP5}, \eqref{estP6}, \eqref{esti4}, \eqref{esti5} and the weak star lower semicontinuity of norms. Therefore, Theorem~\ref{Wellposed} is completely proved.\qed

\section{Continuous dependence and regularity}
\label{Cont}
\setcounter{equation}{0}

In this section, we show report the proofs of the continuous dependence and regularity results.  

\step
Proof of Theorem~\ref{Contdep}

We just have to prove the inequality \eqref{stimaK2} by letting
$(\mu_i,\phi_i,\xi_i)$ be any solution of problem \Pbl\ with the corresponding data $g_i, \mu_{0,i} , \nu_{0,i} ,  \phi_{0,i}$, satisfying \eqref{hyp5}--\eqref{hyp8} for $ i=1,2$. For convenience, within this proof we set 
\Beq
g= g_1 - g_2, \quad \mu_0 = \mu_{0,1}- \mu_{0,2} , \quad \nu_0 = \nu_{0,1}- \nu_{0,2} , \quad 
\phi_0 = \phi_{0,1}- \phi_{0,2} ,
\non
\Eeq  
as well as
\Beq
  \mu = \mu_1-\mu_2 \,, \quad
  \phi= \phi_1-\phi_2 \,, \quad
  \xi = \xi_1-\xi_2.
  \non
\Eeq
Then, taking the differences of the respective equalities~\eqref{prima}--\eqref{seconda} and integrating the one resulting from \eqref{prima} with respect to time, we obtain 
\begin{align}
  &  ( \alpha \dt\mu(t) + \phi(t), v )
  + \iO \nabla (1* \mu) (t) \cdot \nabla v
  = ( \alpha \nu_0 + \phi_0 , v )
  \non
  \\
  & \quad \hbox{\aat\ and every $v\in V$},
  \label{prima-cd}
  \\[1mm]
  &  \tau (\dt \phi(t) , v)
  + \iO \nabla\phi(t) \cdot \nabla v
  + ( \xi(t) , v) = (\mu(t) + g(t) - f'_2 (\phi_1 (t)) +  f'_2 (\phi_2 (t))  , v) 
  \non
  \\
  & \quad \hbox{\aat\ and every $v\in V$},
  \label{seconda-cd}
  \\[2mm]
  &
  \xi_i \in \partial f_1 (\phi_i) , \ i= 1,2, \quad \aeQ ,  
  \label{terza-cd}
  \\[2mm]
  &\mu(0)=\mu_0\,  ,\quad  \phi(0)= \phi_0 \quad \hbox{ a.e.~in }{\Omega} \,.
  \label{cauchy-cd}
\end{align}
Next, we take $v= \mu (t)$ in \eqref{prima-cd} and  $v= \phi (t)$ in \eqref{seconda-cd}, then we add them noting a cancellation of terms and integrate once more with respect to $t$. With the help of~\eqref{cauchy-cd}, we infer that  
\begin{align}
&\frac \alpha 2\| \mu (t)\|_H^2 + \frac 1 2\|\nabla (1* \mu) (t)\|_H^2 +
\frac \tau 2\| \phi(t)\|_H^2 
+ \iint_{Q_t} |\nabla \phi |^2 +
\iint_{Q_t} \xi \phi \non\\
& \leq\,
\frac \alpha 2\|\mu_0\|_H^2 +
\frac \tau 2\|\phi_0 \|_H^2 + 
\iint_{Q_t} ( \alpha \nu_0 + \phi_0) \mu  
\non\\
&\quad{}+ \iint_{Q_t} g \, \phi  - \iint_{Q_t} (f'_2 (\phi_1) -  f'_2 (\phi_2) )\phi
\label{pier6}
\end{align}
for all $t\in (0,T]$. Now, thanks to \eqref{terza-cd} the last term on the \lhs\ is \juerg{nonnegative}. On the \rhs, by the 
Young inequality, we have that  
$$ \iint_{Q_t} ( \alpha \nu_0 + \phi_0) \mu \leq \alpha \int_0^t  \| \mu (s)\|_H^2 ds 
+ C \bigl( \alpha \|\nu_0\|_H^2 +  \alpha^{-1} \|\phi_0\|_H^2 \bigr), $$
and, using also the Lipschitz continuity of \revis{$f'_2$}, it follows that 
$$ \iint_{Q_t} g \, \phi  - \iint_{Q_t} (f'_2 (\phi_1) -  f'_2 (\phi_2) )\phi \leq C \int_0^t  \| \phi (s)\|_H^2 ds + \frac 1 2 \norma{g}_{\L2H}^2.$$
Therefore, we can collect these inequalities and apply the Gronwall lemma to the resultant from \eqref{pier6}
in order to plainly obtain the estimate \eqref{stimaK2}.\qed

\step
Proof of Theorem~\ref{Regularity}

We already know from \eqref{regphi} that $\dt \phi $ is in $\L\infty H \cap \L2V$. Then, in view of 
\eqref{hyp9}, the regularity in \eqref{regmupiu} follows from the variational theory for linear evolution 
problems of second order in time (see, e.g., \cite{DL}). Then, in order to reproduce the estimate in 
\eqref{stimaK3}, let us proceed formally and test equation \eqref{prima} by $-\Delta (\dt \mu (t)) $. By this, we can easily  integrate by parts and also with respect to time. With the help of Young's inequality
we obtain 
\begin{align}
&\frac \alpha 2\iO |\nabla (\dt \mu(t)) |^2 + \frac 1 2 \iO |\Delta \mu(t)|^2
\non\\
&{} = \frac \alpha 2 \iO |\nabla \nu_0|^2 +
\frac 1 2 \iO|\Delta \mu_0 |^2 - \iint_{Q_t} \nabla(\dt\phi) \cdot \nabla (\dt \mu )
\non\\
& \leq\,
\frac \alpha 2\|\nu_0\|_V^2 +
\frac 1 2\|\Delta \mu_0 \|_H^2 
+ \frac1{2\alpha} \iint_{Q_t} |\nabla(\dt\phi)|^2  + \frac \alpha 2 \iint_{Q_t} |\nabla(\dt\mu)|^2,
\label{pier7}
\end{align}
whence the estimate 
\begin{align}
&\alpha^{1/2}\norma{\mu}_{\W{1,\infty}V} 
+ \norma{\mu}_{\L\infty W} 
  \leq C \bigl(1+ \alpha^{- 1/2}\bigr) 
  \label{pier8}
\end{align}
follows from an application of Gronwall's lemma, along with \eqref{hyp9} and \eqref{stimaK1}. 
Having shown \eqref{pier8}, it is now straightforward to compare the terms in \eqref{prima} and to deduce that 
$\alpha \norma{\dtt \phi}_{\L\infty H} $ is bounded by a quantity like the \rhs\ of \eqref{pier8}. 
Thus, we complete the proof of \eqref{stimaK3}.\qed

\step
Proof of Theorem~\ref{Separation}

We first show \eqref{stima-sep}. 
It is already known from \eqref{hyp7} that the initial value of $\phi$, i.e.~$\phiz$, belongs to a compact subset of $D(\partial f_1)=(r_-,r_+) $. 
By the previous proof, we have checked that $\mu$ is bounded in $\L\infty W$,
hence in $\LQ\infty$, as \juerg{follows} from the above estimate \eqref{pier8}. Then, let us rewrite 
equation~\eqref{secondabvp} as 
\Beq
  \tau \dt \phi -\Delta \phi + f_1'(\phi)  = h, \quad \hbox{with }\, h = \mu + g - f_2'(\phi), 
  \quad \aeQ .
  \label{seconda-mp}
\Eeq
The term $\xi$ in \eqref{secondabvp} has been expressed here as $f_1'(\phi)$, as it is allowed by the
assumption~\eqref{hyp11}. Note that the \rhs\ $h$ of \eqref{seconda-mp} is actually bounded in $\LQ\infty$, thanks to \eqref{hyp10} 
and the bound for $\phi$ ensured by \eqref{stimaK1}, along with the Lipschitz continuity of~$f_2'$.

To prove \eqref{stima-sep}, it is enough to derive \juerg{an} $\LQ\infty$-bound for~$f_1'(\phi)$.
Let us outline the argument by proceeding formally {and}
pointing out that \juerg{just} a truncation of \juerg{the} test functions would be needed for a rigorous proof.
We take any $p>2$ and test \eqref{seconda-mp} by $|f_1'(\phi)|^{p-2}f_1'(\phi)$, a function of $\phi$ which is increasing and attains the value 
$0$ at $0$ (cf.~\eqref{hyp2}--\eqref{hyp4}). Then, we integrate from $0$  to $t\in (0,T]$, obtaining
\begin{align}
  &\revis{\tau}\iO \Bigl( \int_0^{\phi(t)} |f_1'(s)|^{p-2}f_1'(s)ds  \Bigr)
  \non 
  \\
  &{}
 + (p-1) \iint_{Q_t} |f_1'(\phi)|^{p-2}f_1''(\phi) |\nabla\phi|^2
  +\iint_{Q_t} |f_1'(\phi)|^p
  \non 
  \\
  &{}= \revis{\tau} \iO \Bigl( \int_0^{\phiz} |f_1'(s)|^{p-2}f_1'(s)ds \Bigr) + 
  \iint_{Q_t} h |f_1'(\phi)|^{p-2}f_1'(\phi) \,.
  \label{pier9}
\end{align}
Note that the first term \revis{on the \lhs\ is nonnegative since $f_1'$ is monotone increasing with $f_1'(0)=0$;
moreover, the second term on the \lhs\ is nonnegative} since the derivative $f_1''$ is 
nonnegative everywhere in $(r_-,r_+) $.  About the \rhs\ we may observe that 
$$ \revis{\tau} \iO \Bigl( \int_0^{\phiz} |f_1'(s)|^{p-2}f_1'(s)ds \Bigr) \leq \revis{\tau} \norma{ f_1'(\phiz)}_\infty^{p-1} 
\norma{ \phiz}_\infty |\Omega|\,, $$ 
and, with \juerg{$p'=p/(p-1)$ and} the help of the Young inequality, that
\begin{align}
  &\iint_{Q_t} h |f_1'(\phi)|^{p-2}f_1'(\phi)
  \leq \norma h_{L^p (Q_t)} \, \norma{\,|f_1'(\phi)|^{p-1}}_{L^{p'} (Q_t)}
  \non
  \\
  &\quad{}
  = \norma h_{L^p (Q_t)} \, \norma{f_1'(\phi)}_{L^p (Q_t)}^{p/p'} 
  \leq \frac 1p \, \norma h_{L^p (Q_t)}^p
  + \frac 1{p'} \, \norma{f_1'(\phi)}_{L^p (Q_t)}^p \,.
  \non
\end{align}
By rearranging \juerg{from} \eqref{pier9}, and taking $t=T$, we infer that 
\begin{align}
  & \norma{f_1'(\phi)}_{L^p (Q)} \leq  \Bigl( p \revis{\tau} \norma{ f_1'(\phiz)}_\infty^{p-1} 
\norma{ \phiz}_\infty |\Omega| +  \norma h_{L^p (Q)}^p \Bigr)^{1/p}  
  \non
  \\
  &\quad{}
  \leq  \bigl( p \revis{\tau} \norma{ f_1'(\phiz)}_\infty^{p-1} 
\norma{ \phiz}_\infty |\Omega| \bigr)^{1/p}  +  \norma h_{L^p (Q)}.
\non
\end{align}
Then, \revis{passing to the limit as $p\to +\infty$ in the above chain of inequalities}, we conclude~that 
$$\norma{f_1'(\phi)}_{\juerg{L^\infty} (Q)} 
\leq  \norma{ f_1'(\phiz)}_\infty
  +  \norma h_{\juerg{L^\infty} (Q)}, $$
which ensures the validity of \eqref{stima-sep}, for some constants  $r_*, r^* $ as in the statement.

Next, we argue in order to prove the continuous dependence estimate \eqref{stimaK4}. We use the same notation as in the proof of Theorem~\ref{Contdep}, so that 
\Beq
g= g_1 - g_2, \quad \mu_0 = \mu_{0,1}- \mu_{0,2} , \quad \nu_0 = \nu_{0,1}- \nu_{0,2} , \quad 
\phi_0 = \phi_{0,1}- \phi_{0,2}, 
\non
\Eeq  
and 
\Beq
  \mu = \mu_1-\mu_2 \,, \quad
  \phi= \phi_1-\phi_2 \,,
  \non
\Eeq
where $(\mu_i,\phi_i,\xi_i )$, with $\xi_i = f_1'(\phi_i)$, is the solution to problem \Pbl\ corresponding to  $g_i, \mu_{0,i} , \nu_{0,i} ,  \phi_{0,i}$, $i=1,2$, these data satisfying \eqref{hyp5}--\eqref{hyp8}.
Then, taking the differences of the respective equalities~\eqref{prima}--\eqref{seconda}, we obtain 
\begin{align}
  & \alpha \langle \dtt\mu(t) , v \rangle  
  + \iO \nabla\mu(t) \cdot \nabla v
  = - ( \dt\phi(t), v )
  \non
  \\
  & \quad \hbox{\aat\ and every $v\in V$},
  \label{prima-cd2}
  \\[1mm]
  &  \tau (\dt \phi(t) , v)
  + \iO \nabla\phi(t) \cdot \nabla v
  = (\mu(t) + g(t) - f' (\phi_1 (t)) +  f' (\phi_2 (t))  , v) 
  \non
  \\
  & \quad \hbox{\aat\ and every $v\in V$},
  \label{seconda-cd2}
  \\[2mm]
  &\mu(0)=\mu_0\, ,\quad  (\dt \mu) (0)=\nu_0 \, ,\quad  \phi(0)= \phi_0 \quad \hbox{ a.e.~in }{\Omega} \,,
  \label{cauchy-cd2}
\end{align}
where we have used $f' = f_1' + f_2' $ in equation~\eqref{seconda-cd2} noting that the estimate \eqref{stima-sep} and 
the assumptions~\eqref{hyp11} on~$f_1$
now allow us to take $f'$ as a global Lipschitz continuous function
(with Lipschitz constant depending on $\alpha$). In view of \eqref{stimaK2}, 
\juerg{we have} for the term on the \rhs\ of~\eqref{seconda-cd2} that 
\begin{align}
&\norma{\mu + g - f' (\phi_1) +  f' (\phi_2)}_{\L2 H}
\leq \norma{\mu}_{\L2 H} + 
\norma{g}_{\L2 H} + 
C_\alpha \norma{\phi}_{\L2 H}
\non \\
&\quad {}\leq C_\alpha \Bigl( \norma{g}_{\L2H} + \norma{\mu_{0}}_{H} +
\norma{\nu_{0}}_{H}  + \norma{\phi_0 }_{H}   \Bigr),
  \label{pier10}
\end{align}
where the constants are denoted by $C_\alpha$ since they depend on $\alpha$ as well.
Then, using the standard parabolic regularity estimate (see, e.g., \cite{Lions} or \cite{DL}) 
\begin{align}
& \norma{\phi}_{\H1H\cap\L\infty V\L2W} 
\non\\
&\leq C \bigl(\norma{\mu + g - f' (\phi_1) +  f' (\phi_2)}_{\L2 H} + 
 \norma{\phi_0 }_{V}\bigr)\non 
 \end{align}
for \juerg{the} $\phi$ solution \juerg{to} \eqref{seconda-cd2} with the respective initial condition, 
it is straightforward to deduce that 
\begin{align}
&\norma{\phi}_{\H1H\cap\L\infty V\L2W} 
\non
\\
&{}\leq C_\alpha \bigl( \norma{g}_{\L2H} + \norma{\mu_{0}}_{H} +
\norma{\nu_{0}}_{H}  + \norma{\phi_0 }_{V}   \bigr).
\label{pier11}
\end{align} 
Next, we can choose $v= \dt \mu $ in \eqref{prima-cd2}, integrate with respect to time, and infer that
\begin{align}
&\frac \alpha 2\| \dt \mu (t)\|_H^2 + \frac 1 2\|\nabla  \mu (t)\|_H^2 
\non
\\
&{} \leq\,
\frac \alpha 2\|\nu_0\|_H^2 +
\frac 1 2\|\nabla \mu_0 \|_H + 
\int_0^t \| \dt \phi (s)\|_H \| \dt \mu (s)\|\juerg{_H} ds .
\label{pier12}
\end{align}
Then, \juerg{first applying Young's inequality to the last term and then Gronwall's lemma}, we arrive at the estimate
\begin{align}
& \norma{\mu}_{\W{1,\infty}H\cap\L\infty V} 
\leq C_\alpha  \bigl( \norma{\mu_{0}}_{V} + \norma{\nu_{0}}_{H} + \norma{\dt \phi }_{\L2 H} \bigr), 
\non 
\end{align}
whence from \eqref{pier11} it is clear that  
\begin{align}
& \norma{\mu}_{\W{1,\infty}H\cap\L\infty V} 
\leq C_\alpha \bigl( \norma{g}_{\L2H} + \norma{\mu_{0}}_{V} +
\norma{\nu_{0}}_{H}  + \norma{\phi_0 }_{V}   \bigr).
\non 
\end{align}
By this estimate, comparison of the terms in \eqref{prima-cd2} yields 
\begin{align}
&\alpha \norma{\dtt \mu }_{\L2 {V^*}} \leq \norma{\nabla \mu}_{\L2 H } + \norma{\dt \phi }_{\L2 H}
\non
\\
&\quad {}
\leq C_\alpha \bigl( \norma{g}_{\L2H} + \norma{\mu_{0}}_{V} +
\norma{\nu_{0}}_{H}  + \norma{\phi_0 }_{V}   \bigr),
  \label{pier13}
\end{align}
so that \eqref{stimaK4} is completely proved.\qed

\section{Asymptotic analysis}
\label{Asy}
\setcounter{equation}{0}

This section is devoted to the study of the asymptotic behavior of the problem \State\ as $\alpha $ approaches $0$. We allow the initial data for $\mu$ and $\dt \mu$, as well as the source term $g$, to depend on $\alpha$, while we keep fixed $\phiz$, the initial value of $\phi$, for reasons of simplicity in front of
restrictions like \eqref{hyp7} and \eqref{hyp8} for $\phi_0$. 

Thus, for $0< \alpha \leq 1$, we consider families of data $g_\alpha, \mu_{0,\alpha} ,  \nu_{0,\alpha}$ such that 
\begin{align}
& \{ g_\alpha \} \, \hbox{ is uniformly bounded in }\,  \H1H
\non \\
& \qquad \hbox{and strongly converges to $g$ in } \, \L2H \, \hbox{ as }
\alpha\searrow 0, \label{hyp13}\\[1mm]
&\{ \mu_{0,\alpha}\} \, \hbox{ is uniformly bounded in }\,  V ,\label{hyp14}\\[1mm]
& \{ \nu_{0,\alpha}\} \, \hbox{ is uniformly bounded in }\,  H . \label{hyp15}
\end{align} 
Of course, it follows from \eqref{hyp13} that $g\in \H1H$ and $\dt g_\alpha \to \dt g$ weakly in $\L2H$. 
We can state the following convergence result. 

\Bthm
\label{Convergence}
Assume that \eqref{hyp1}--\eqref{hyp4}, \eqref{hyp7}--\eqref{hyp8}, \eqref{hyp13}--\eqref{hyp15} are fulfilled.
For all $\alpha \in (0,1]$, let the triple $(\mu_\alpha , \phi_\alpha , \xi_\alpha)$, with 
\begin{align}
  & \mu_\alpha \in \W{2,\infty}{V^*}\cap \W{1,\infty}H \cap \L\infty V,
  \label{regmual}
  \\
  & \phi_\alpha \in \W{1,\infty}H \cap \H1V \cap \L\infty W,
  \label{regphial}
  \\
  & \xi_\alpha \in \L\infty H,
  \label{regxial}
\end{align}
be the solution \juerg{to} the initial value problem
\begin{align}
  & \alpha \< \dtt\mua(t) , v > 
  + ( \dt\phia(t), v )
  + \iO \nabla\mua(t) \cdot \nabla v
  = 0
  \non
  \\
  & \quad \hbox{\aat\ and every $v\in V$},
  \label{primaal}
  \\[1mm]
  &  \tau (\dt \phia(t) , v)
  + \iO \nabla\phia(t) \cdot \nabla v
  + ( \xia(t) +  f'_2 (\phia(t)), v) = (\mua(t) + g_\alpha (t), v) 
  \non
  \\
  & \quad \hbox{\aat\ and every $v\in V$},
  \label{secondaal}
  \\[2mm]
  &
  \xia \in \partial f_1 (\phia) \quad \aeQ ,  
  \label{terzaal}
  \\[2mm]
  &\mua(0)=\muza \, ,\quad  (\dt \mua) (0)=\nuza \, ,\quad  \phia(0)= \phiz \quad \hbox{ a.e.~in }{\Omega} \,.
  \label{cauchyal}
\end{align}
Then there exists a triple $(\mu,\phi,\xi)$ such that,
for some subsequence $\alpha_k$ tending to $0$, there holds  
\begin{align}
\label{conv7}
&\mu_{\alpha_k} \to\mu \quad \mbox{weakly star in }\,\L\infty V,\\
\label{conv8}
&\alpha_k\mu_{\alpha_k} \to 0 \quad \mbox{weakly star in }\,\W{2,\infty}{V^*} \,  
\mbox{ and strongly in }\,\W{1,\infty}H,\\
\label{conv9}
&\phi_{\alpha_k}\to\phi \quad \mbox{weakly star in } \,W^{1,\infty}(0,T;H)\cap H^1(0,T;V)\cap L^\infty(0,T;W)
\non
\\
&\qquad 
\mbox{ and strongly in }\,\C0V \cap C^0(\overline Q),\\
\label{conv10}
&\xi_{\alpha_k}\to\xi \quad \mbox{weakly star in } L^\infty(0,T;H).
\end{align}
Moreover, $(\mu,\phi,\xi)$ is a solution \juerg{to} the viscous Cahn--Hilliard system 
\begin{align}
  &  ( \dt\phi(t), v )
  + \iO \nabla\mu(t) \cdot \nabla v
  = 0
  \non
  \\
  & \quad \hbox{\aat\ and every $v\in V$},
  \label{primaCH}
  \\[1mm]
  &  \tau (\dt \phi(t) , v)
  + \iO \nabla\phi(t) \cdot \nabla v
  + ( \xi(t) +  f'_2 (\phi(t)), v) = (\mu(t) + g(t), v) 
  \non
  \\
  & \quad \hbox{\aat\ and every $v\in V$},
  \label{secondaCH}
  \\[2mm]
  &
  \xi \in \partial f_1 (\phi) \quad \aeQ ,  
  \label{terzaCH}
  \\[2mm]
  & \phi(0)= \phi_0 \quad \hbox{ a.e.~in }{\Omega} \,.
  \label{cauchyCH}
\end{align}
\Ethm

\Bdim
A closer inspection of the proof of Theorem~\ref{Wellposed} in Section~\ref{Exi} reveals that the estimate~\eqref{stimaK1} still holds under the conditions~\eqref{hyp13}--\eqref{hyp15}. Then, by a 
standard weak star compactness argument, we deduce the existence of a subsequence $\alpha_k\searrow 0$ 
and a triple $(\mu,\phi,\xi)$ such that \eqref{conv7}--\eqref{conv10} hold. In fact, the strong convergence  
property in \eqref{conv9} is a consequence of the compactness result reported in \cite[Sect.~8,~Cor.~4]{Simon}. Moreover, by the Lipschitz continuity of $f_2'$, we also have that 
\begin{align}
\label{conv11}
&f_2' (\phi_{\alpha_k}) \to f_2'(\phi) \quad \mbox{strongly in } \,C^0(\overline Q)\cap \C0H\, .
\end{align}
Then, we can pass to the limit in \eqref{primaal}, \eqref{secondaal}, and the third condition 
in~\eqref{cauchyal}, all written for $\alpha_k$, and easily obtain \eqref{primaCH}, \eqref{secondaCH}, \eqref{cauchyCH}. Recovering \eqref{terzaCH} from \eqref{terzaal} is straightforward, due to the weak convergence of $\xi_{\alpha_k}$ and the strong convergence of $\phi_{\alpha_k}$ in $\L2H$, along with the maximal monotonicity of $\partial f_1$ (see, e.g.,~\cite[Cor.~2.4, p.~41]{Barbu}). \juerg{This concludes} the proof.
\Edim

\Brem
Note that Theorem~\ref{Convergence} implicitly yields an existence result for solutions to the viscous 
Cahn--Hilliard system~\eqref{primaCH}--\eqref{cauchyCH}. The found solution is already a regular and strong solution\revis{. Indeed,} the component 
$\phi$ is in $W^{1,\infty}(0,T;H)\cap H^1(0,T;V)\cap L^\infty(0,T;W)$ and therefore also in  
\revis{$\C0{H^{2-\delta} (\Omega)}$ for all $\delta>0$, by compactness (see, e.g., \cite[Sect.~8,~Cor.~4]{Simon}).
As a conseqence of the dimension $N\leq3$ of the space and of the Sobolev embedding theorems, it turns out that $\phi \in C^0(\overline Q)$. Moreover,} from \eqref{conv7}, 
a comparison in \eqref{primaCH}, 
and the elliptic regularity theory, it turns out that $\mu \in \L\infty W \cap \L2{H^3(\Omega)}$ in addition, 
\juerg{so that $\mu \in \revis{\L\infty{C^0(\overline{\Omega})}}$,} in particular. Both \juerg{the} equalities \eqref{primaCH} and \eqref{secondaCH} can be equivalently rewritten as the equations

\begin{align}
   &\dt \phi - \Delta \mu = 0  \quad \aeQ ,
   \label{pier14}
   \\
  & \tau \dt \phi -\Delta \phi + \xi + f_2'(\phi) = \mu + g 
  \quad \aeQ, 
  \label{pier15}
\end{align}
plus the homogeneous boundary conditions
\Beq
\dn \mu = \dn\phi = 0 \quad \mbox{a.e. on }\,\Sigma.
\label{pier16}
\Eeq
The mentioned regularity for  $(\mu,\phi,\xi)$ is exactly the same as in \cite[Thm.~2.2]{CoSpTr}, where \juerg{a} slightly more general system is investigated. However, the existence of a less regular solution can also be proved, along with the uniqueness of the component $\phi$ of the solution, as it results for instance from 
\cite[Thm.~2.5]{CGM}. Please note that in general uniqueness cannot be expected for $\xi$ and $\mu$ unless $\partial f_1$ is single-valued (like e.g. the case considered in \eqref{hyp11}); otherwise, only the difference $\xi - \mu $ is uniquely determined from~\eqref{pier15}. 

\Erem 

\Brem
By the uniqueness property for the component $\phi$, which is pointed out in the previous remark, we infer 
that not only a subsequence $\{\phi_{\alpha_k} \}$ but the entire family $\{\phia\}_{\alpha\in (0,1]}$ 
converges to $\phi$ in the sense of \eqref{conv9} as $\alpha \searrow 0$. 
\Erem

The next result is devoted to an error estimate of the difference $\phia - \phi$ in certain norms and in terms of the parameter $\alpha$.  

\Bthm
\label{Errest}
Under the same assumptions as in Theorem~\ref{Convergence}, we let $(\mu_\alpha , \phi_\alpha , \xi_\alpha)$ denote the solution to \eqref{primaal}--\eqref{cauchyal}, for $\alpha \in (0,1]$, and $(\mu,\phi,\xi)$ be the solution to \eqref{primaCH}--\eqref{cauchyCH} found by the asymptotic limit in \eqref{conv7}--\eqref{conv10}. Then there is a constant $K_5>0$, which depends on the structure of the system but is independent of $\alpha$, such that
\begin{align}
&\alpha^{1/2} \norma{\mu_a}_{\L\infty H}
+ \norma{\nabla (1* (\mua-\mu)) }_{\L\infty H}
+ \norma{\phia-\phi}_{\L\infty H\cap\L2V}
\non \\
&{}\leq K_5 \, \Bigl(\alpha^{1/4} + \norma{g_\alpha -g}_{\L2H} \Bigr).
  \label{stimaK5}
\end{align} 
\Ethm

\Bdim
We argue similarly as in the proof of Theorem~\ref{Contdep}. We take the difference of \eqref{primaal} and \eqref{primaCH}, then we integrate with respect to time with the help of \eqref{cauchyal} and \eqref{cauchyCH}. We obtain
\begin{align}
  &  ( (\phia - \phi)(t), v )
  + \iO \nabla (1*(\mua- \mu)) (t) \cdot \nabla v
  = ( \alpha \, \nuza -  \alpha \, \dt \mua(t), v )
  \non
  \\
  & \quad \hbox{\aat\ and every $v\in V$}.
  \label{pier17}
\end{align}
At the same time, we subtract \eqref{secondaCH} from \eqref{secondaal} and have that 
\begin{align}
  &  \tau (\dt (\phia - \phi) (t) , v)
  + \iO \nabla( \phia - \phi)(t) \cdot \nabla v
  + ( (\xia - \xi) (t) , v) 
  \non
  \\
  &{}
  = ((\mua - \mu) (t) + (g_\alpha - g) (t) - f'_2 (\phia (t)) +  f'_2 (\phi (t))  , v) 
  \non
  \\
  & \quad{} \hbox{\aat\ and every $v\in V$}\,.
  \label{pier18}
\end{align}
Then we take $v= (\mua-\mu) (t) $ in \eqref{pier17} and  $v= (\phia-\phi) (t) $ in \eqref{pier18}, sum up  noting that a cancellation occurs, and integrate with respect to $t$. Since the product $ (\xia - \xi) (\phia-\phi)$ is nonnegative due to \eqref{terzaal}, \eqref{terzaCH} and the monotonicity of $\partial f_1$, we easily derive the inequality 
\begin{align}
&\frac \alpha 2\| \mua (t)\|_H^2 + \frac 1 2\|\nabla  (1*(\mua- \mu)) (t)\|_H^2 +
\frac \tau 2\| (\phia-\phi)(t)\|_H^2 
+ \iint_{Q_t} |\nabla (\phia -\phi) |^2 
\non\\
& \leq\,
\frac \alpha 2\|\muza \|_H^2 +
\iint_{Q_t}  \alpha \, \nuza  (\mua- \mu) + \iint_{Q_t} \alpha \, \dt \mua \, \mu 
\non\\
&\quad{}+ \iint_{Q_t} (g_\alpha - g)  (\phia -\phi)  - \iint_{Q_t} (f'_2 (\phia) -  f'_2 (\phi) ) (\phia -\phi)
\label{pier19}
\end{align}
for all $t\in (0,T]$. Now, we recall the boundedness properties \eqref{hyp14} and \eqref{hyp15}, the estimate  \eqref{stimaK1} for $\norma{\mua}_{\L\infty H}$ and 
$\alpha^{1/2} \norma{\dt\mua}_{\L\infty H} $, \juerg{as well as} the regularity  $\mu \in \L\infty H$, in order to deduce that
\begin{align}
&\frac \alpha 2\|\muza \|_H^2 +
\iint_{Q_t}  \alpha \, \nuza  (\mua- \mu) + \iint_{Q_t} \alpha \, \dt \mua \, \mu 
\non\\
&{}\leq 
\alpha\norma{\muza}_{H}^2  + \revis{C}
\alpha \norma{\nuza}_{H} + 
C \alpha^{1/2} \leq C \alpha^{1/2} .
\non 
\end{align}
In addition, by virtue of the Lipschitz continuity of $f_2'$ and the Young inequality, we have that 
\begin{align}
&\iint_{Q_t} (g_\alpha - g)  (\phia -\phi)  - \iint_{Q_t} (f'_2 (\phia) -  f'_2 (\phi) ) (\phia -\phi)
\non
\\
&{}\leq \norma{g_\alpha -g}_{\L2H}^2 + C \int_0^t \| (\phia -\phi)(s)\|^2 ds.
\non
\end{align}
Then, collecting the above computations in \eqref{pier19} and applying Gronwall's lemma, the estimate \eqref{stimaK5} follows. 
\Edim

We finally notice that \eqref{stimaK5} gives an error estimate, in particular, for 
$$ \norma{\phia-\phi}_{\L\infty H\cap\L2V}$$
of order $1/4$, provided that the convergence
of $\norma{g_\alpha -g}_{\L2H}$ to $0$ is at least of this order. 
 

\smallskip

\section*{Acknowledgments}
This research has been supported by the MIUR-PRIN Grant
2020F3NCPX ``Mathematics for industry 4.0 (Math4I4)''. 
In addition, PC indicates his affiliation 
to the GNAMPA (Gruppo Nazionale per l'Analisi Matematica, 
la Probabilit\`a e le loro Applicazioni) of INdAM (Isti\-tuto 
Nazionale di Alta Matematica). During the preparation of this paper, individual visits to the WIAS in Berlin
and to the Dipartimento di Matematica ``F. Casorati'' in Pavia have been exchanged by the two authors: the warm hospitality of the two institutions is gratefully acknowledged. 
\bigskip


\End{document}